\font\Bbb=msbm10 
\def\matC{\hbox{\Bbb C}}

\font\Gros=cmbx10 scaled\magstep1

\magnification \magstep1 

\null  
\bigskip
\bigskip
\centerline{\Gros La pseudodistance de Carath\'eodory  } 
\smallskip
\centerline{\Gros sur des ouverts embo\^{\i}t\'es} 
\smallskip

\medskip
\centerline{\Gros Jean-Pierre Vigu\'e}
\bigskip

{\bf1. Introduction} 

\medskip
On consid\`ere une vari\'et\'e analytique complexe $ X $ de dimension finie 
et un ouvert $ U $ de $ X$. Pour les pseudodistances de Carath\'eodory $c_{X}$ et 
$c_{U}$,
 on sait que, $ \forall x,y \in U,$ 
$$c_{X}(x,y)\leq c_{U}(x,y).$$
 Cependant, sous des hypoth\`eses un peu plus fortes, nous montrerons
 qu'il existe une constante $ k<1 $ telle que, $ \forall x,y \in U,$ 
$$c_{X}(x,y)\leq kc_{U}(x,y).$$
 En particulier, si $ X $ est $ c_{X}$-hyperbolique, ce qui signifie
 que $ c_{X} $ est une distance qui d\'efinit la topologie de $ X$, ceci
 entra\^{\i}ne que toute application holomorphe $ f:X\longrightarrow
 X $ telle que $ f(X)\subset U $ admet un point fixe unique. 

Ce m\^eme probl\`eme pour la pseudodistance int\'egr\'ee de Carath\'eodory
 a d\'ej\`a \'et\'e \'etudi\'e par H.-J. Reiffen [5], C. Earle and R.
 Hamilton [1] et J.-P. Vigu\'e [6]. 

\medskip

{\bf2. Rappel sur les pseudodistances invariantes} 

\medskip

Sur le disque unit\'e $ \Delta\subset{\matC}$, on d\'efinit la distance
 de Poincar\'e $ \omega $ par la formule : 
pour tous $ z $ et $ w $ appartenant \`a $ \Delta$,   
$$\omega(z,w)= \hbox{\rm\ tanh }^{-1}\Big|\,{z-w      \over 1-\overline{w}z
     }\,\Big|\hbox{\rm.}$$
Ici, ${\rm tanh }^{-1}$ d\'esigne l'inverse de la fonction tangente hyperbolique.

On d\'efinit alors la pseudodistance de Carath\'eodory $ c_{X} $ sur
 une vari\'et\'e analytique complexe $ X $ par la formule : pour tous
 $ x $ et $ y $ appartenant \`a $ X$, 

$$c_{X}(x,y)=\hbox{\rm sup}_{\varphi\in H(X,\Delta)} \omega(\varphi(x),\varphi(y)),$$
 o\`u $ H(X,\Delta) $ d\'esigne l'ensemble des applications holomorphes de $ X $ dans
 le disque-unit\'e $ \Delta$. 
 Comme $ \Delta $ est homog\`ene et
 que les automorphismes analytiques de $ \Delta $ sout des isom\'etries
 pour $ \omega$, on peut supposer de plus que $ \varphi(x)=0$. 

On v\'erifie facilement que $c_{\Delta}=\omega$ 
et que $ c_{X} $ est une pseudodistance invariante,
 ce qui signifie que toute application holomorphe $ f:X\longrightarrow
 Y $ v\'erifie, pour tous $ x $ et $ y $ appartenant \`a $ X$, 
$$c_{Y}(f(x),f(y))\leq c_{X}(x,y)$$
 et est donc contractante (au sens large)(voir [2,3 et 4]). 
\medskip
{\bf3. Lemme de convexit\'e} 
\medskip
On a le lemme de convexit\'e suivant. 

{\bf Lemme 3.1. }{\it Soit \/}$ r<1 $ {\it un nombre positif. Alors,
 pour tout \/}$ x\geq0${\it,\/}  

$${\rm tanh}^{-1}(r{\rm tanh }(x))\leq rx,$$ 

\noindent {\it o\`u ${\rm  tanh}$ d\'esigne la fonction 
tangente hyperbolique et ${\rm tanh}^{-1}$
  {\it son inverse.\/} 

{\it D\'emonstration. \/} \rm On a : 

$$({\rm tanh}^{-1})'(x)=1/(1-x^{2}), $$
 et 

$$({\rm tanh}^{-1})''(x)=2x/(1-x^{2})^{2}.$$ 

Par suite,  ${({\rm tanh}^{-1})}''(x) $ est positif pour $ x\geq0$, et la fonction
 $ {\rm tanh}^{-1} $ est convexe sur $ [0,1[$. On peut alors utiliser l'in\'egalit\'e
 de convexit\'e au point $ r {\rm tanh }(x) $ sur le segment $ [0,{\rm tanh }x]$, et
 on trouve 

${\rm tanh}^{-1}(r {\rm tanh }(x))\leq r{\rm tanh }^{-1}({\rm tanh }(x))\leq rx$, 

\noindent et le r\'esultat est d\'emontr\'e. 
\medskip

{\bf4. La pseudodistance de Carath\'eodory sur des ouverts embo\^{\i}t\'es}
  \medskip

Nous pouvons montrer le th\'eor\`eme suivant. 

{\bf Th\'eor\`eme 4.1. }{\it Soit \/}$ X $ {\it une vari\'et\'e analytique
 complexe. Soit \/}$ U $ {\it un ouvert de \/}$ X $ {\it born\'e pour la
 pseudodistance de Carath\'eodory \/}$ c_{X}${\it. Soit \/}$ M=sup_{x\in
 U,y\in U}c_{X}(x,y)${\it et soit \/}$ k={\rm tanh }M<1$\it. Alors, pour tout
 \/}$ x\in U${\it, pour tout \/}$ y\in U${\it,\/} 

{\it c\/}$_{\hbox{\it X\/}}${\it(x,y)\/}$\leq${\it kc\/}$_{\hbox{\it
 U\/}}${\it(x,y).\/} 
\medskip
{\it D\'emonstration.\/} Remarquons que $ k=0 $ signifie que $ \forall
 x\in U$, $ \forall y\in U$, $ c_{X}(x,y)=0$, et le r\'esultat est \'evident.
 On peut donc supposer $ k>0 $ et soit $ x $ un point de $ U$.
 Soit $ f:X\longrightarrow\Delta $ une application holomorphe
 telle que $ f(x)=0$. Du fait que $ f $ est contractante pour la pseudodistance
 de Carath\'eodory, on d\'eduit que, $ \forall y\in U$, 
$$c_{\Delta}(f(x),f(y))\leq c_{X}(x,y)\leq M.$$
 Par suite, $ |f(y)|\leq {\rm tanh }M=k<1$. 

\noindent On en d\'eduit que la fonction $ \varphi $ d\'efinie par $ \varphi(z)=(1/k)f(z)
 $ envoie $ U $ dans le disque-unit\'e $ \Delta$. Ceci implique que 

$$c_{U}(x,y)\geq c_{\Delta}((1/k)|f(x)|,(1/k)|f(y)|),$$ 

$$c_{U}(x,y)\geq {\rm tanh }^{-1}((1/k)|f(y)|).$$ 

\noindent En prenant la borne sup\'erieure pour toutes les applications holomorphes
 $ f:X\longrightarrow\Delta $ telles que $ f(x)=0$, on trouve

$$c_{U}(x,y)\geq {\rm tanh} ^{-1}((1/k){\rm tanh }c_{X}(x,y)).$$

\noindent Par suite, 

$${\rm tanh }c_{U}(x,y)\geq(1/k){\rm tanh }c_{X}(x,y) $$ 

\noindent et 

$${\rm tanh}^{-1}(k{\rm tanh }c_{U}(x,y))\geq c_{X}(x,y).$$ 

\noindent D'apr\`es le lemme 3.1, on a 
$kc_{U}(x,y)\geq c_{X}(x,y)$, et le th\'eor\`eme est d\'emontr\'e.

En particulier, si $ U $ est relativement compact dans $ X$, alors $ U
 $ est born\'e pour $ c_{X}$. On peut donc utiliser le th\'eor\`eme pr\'ec\'edent
 et on a le corollaire suivant. 

{\bf Th\'eor\`eme 4.2. }{\it Soit \/}$ X $ {\it une vari\'et\'e analytique
 complexe. Soit \/}$ U $ {\it un ouvert de \/}$ X $ {\it relativement compact
 dans \/}$ X${\it, soit \/}$ M=sup_{x\in U,y\in U}c_{X}(x,y)${\it et
 soit \/}$ k={\rm tanh }M<1${\it. Alors, pour tout \/}$ x\in U${\it, pour tout
 \/}$ y\in U${\it,\/} 

{\it c\/}$_{\hbox{\it X\/}}${\it(x,y)\/}$\leq${\it kc\/}$_{\hbox{\it
 U\/}}${\it(x,y).\/}    

\medskip

On en d\'eduit le corollaire  4.3.

{\bf Corollaire 4.3.} \it Soit $ X $ une vari\'et\'e $c_X $-hyperbolique
 et soit $ f:X\longrightarrow
 X $ une application holomorphe telle que $ f(X) $ soit relativement compact
 dans $ X$. Alors la suite des it\'er\'ees $ f^{n} $ converge, pour la
 topologie compacte ouverte vers une application constante $ z\mapsto b$,
 o\`u $ b $ est l'unique point fixe de $ f$.

{\it D\'emonstration.}  \rm
On d\'eduit facilement du fait que $ f(X) $ est relativement compact
 dans $ X$ qu'il existe un ouvert $ U $ relativement compact dans $ X $ tel
 que $ f(X)\subset U\subset X$. Comme $ f $ envoie $ X $ dans $ U$, on
 a, pour tout $ x\in X$, pour tout $ y\in X$, 
$$c_{U}(f(x),f(y))\leq c_{X}(x,y).$$
 D'apr\`es le th\'eor\`eme 4.2, il existe une constante $ k<1 $ telle
 que 
$$c_{X}(f(x),f(y))\leq kc_{U}(f(x),f(y)).$$

\noindent  On trouve alors 
$$c_{X}(f(x),f(y))\leq kc_{X}(x,y).$$

\noindent  En utilisant ce r\'esultat pour $ y=f(x)$, et en it\'erant, on en
 d\'eduit que, pour tout $ n>0$, 
$$c_X(f^{n}(x),f^{n+1}(x))\leq k^{n}c_X(x,f(x)).$$
 De fa\c con tout \`a fait classique, on en d\'eduit que, pour tout
 $ a\in X$, la suite des it\'er\'ees $ (f^{n}(a)) $ est une suite de Cauchy
 pour $c_X$ sur $ \overline{f(X)} $ qui est compact. Elle converge donc
 vers $ b $ qui est l'unique point fixe de $ f$. 

\medskip

{\bf Bibliographie}
\medskip

1. C. Earle and R. Hamilton. A fixed point theorem for holomorphic mappings. 
Proc. Symposia Pure Math., {\bf16} (1970),
 61--65.
\medskip

2. T. Franzoni and E. Vesentini. Holomorphic maps and invariant distances.
 Notas de Matematica [Mathematical Notes], {\bf69}. North-Holland Publishing
 Co, Amsterdam, 1980. 
\medskip

3. M. Jarnicki and P. Pflug. Invariant distances and metrics in complex
 analysis. de Gruyter Expositions in Mathematics, {\bf9}, Walter de
 Gruyter Co, Berlin, 1993. 
\medskip

4. S. Kobayashi. Hyperbolic complex spaces, Grundlehren 
der Mathematischen Wissenschaften [Fundamental Principles 
of Mathematical Sciences], {\bf 318}, Springer-Verlag, Berlin, 1998.
\medskip

5. H.-J. Reiffen. Die Carath\'eodorysche Distanz und ihre zugeh\"orige
 Differentialmetrik. Math. Ann.  {\bf161}  (1965), 315--324.
\medskip

6. J.-P. Vigu\'e. Distances invariantes et points fixes d'applications holomorphes.
 Bull. Sci. Math. {\bf136} (2012), 12--18.

\medskip
\bigskip 

Jean-Pierre Vigu\'e

Universit\'e de Poitiers 

Math\'ematiques 

SP2MI, BP 30179
 
86962 FUTUROSCOPE  

e-mail : vigue@math.univ-poitiers.fr

ou jp.vigue@orange.fr

\bye